\documentclass[12pt]{amsart}
\usepackage[mathscr]{eucal}
\usepackage{amscd}
\usepackage{verbatim}

\def\Bbb{\mathbb}
\def\frak{\mathfrak}

\newenvironment{pf*}[1]{\proof[#1]}{\endproof}
\newtheorem{Theorem}[equation]{Theorem}
\newtheorem{Corollary}[equation]{Corollary}
\newtheorem{Lemma}[equation]{Lemma}

\theoremstyle{definition}
\newtheorem{Definition}[equation]{Definition}
\newtheorem{Example}[equation]{Example}

\theoremstyle{remark}
\newtheorem{Remark}[equation]{Remark}

\numberwithin{equation}{section}

\newcommand{\thmref}[1]{Theorem~\ref{#1}}
\newcommand{\secref}[1]{\S\ref{#1}}
\newcommand{\lemref}[1]{Lemma~\ref{#1}}

\newcommand{\corref}[1]{Corollary~\ref{#1}}

\newcommand{\defeq}{\overset{\operatorname{\scriptstyle def.}}{=}}
\newcommand{\C}{{\Bbb C}}
\newcommand{\Z}{{\Bbb Z}}
\newcommand{\Q}{{\Bbb Q}}

\newcommand{\g}{{\frak g}}

\newcommand{\ve}{\varepsilon}

\newcommand{\Uq}{{\mathbf U}_q(\mathfrak g)} 

\newcommand{\Ua}{{\mathbf U}_q(\widehat{\mathfrak g})} 
\newcommand{\Ul}{{\mathbf U}_q({\mathbf L}{\mathfrak g})} 

\newcommand{\Lg}{\mathbf L\g}

\newcommand{\Y}[2]{Y_{#1,#2}}
\newcommand{\Yq}[2]{
\ifx#21 \Y{#1}{aq}
\else\ifx#20 \Y{#1}{a}
\else \Y{#1}{aq^{#2}}
\fi\fi}

\newcommand{\bfR}{\mathbf R}

\newcommand{\qch}{\chi_{q}}
\newcommand{\qtch}{\chi_{q,t}}
\newcommand{\tilqtch}{\widetilde{\chi_{q,t}}}

\newcommand{\PKr}[3]{
\ifx#31 P_{#2,aq}^{(#1)}
\else\ifx#30 P_{#2,a}^{(#1)}
\else P_{#2,aq^{#3}}^{(#1)}
\fi\fi}

\newcommand{\Kr}[3]{
\ifx#30 W_{#2,a}^{(#1)}
\else\ifx#31 W_{#2,aq}^{(#1)}
\else W_{#2,aq^{#3}}^{(#1)}
\fi\fi}

\newcommand{\QK}[2]{Q^{(#1)}_{#2}}
\newcommand{\maxexp}{\operatorname{r}}

\newcommand{\yl}{14pt}

\newcommand{\twoY}[2]{
\renewcommand{\arraycolsep}{0pt}
\begin{array}{|c|}
   \hline \hbox to \yl{\hfill$\scriptstyle #1$\hfill} \\
   \hline \hbox to \yl{\hfill$\scriptstyle #2$\hfill} \\
   \hline
\end{array}
}

\begin{document}
\title[Kirillov-Reshetikhin modules of quantm affine algebras]
{$t$--analogs of $q$--characters of Kirillov-Reshetikhin modules of
quantm affine algebras
}
\author{Hiraku Nakajima}
\address{Department of Mathematics, Kyoto University, Kyoto 606-8502,
Japan
}
\email{nakajima@kusm.kyoto-u.ac.jp}
\urladdr{http://www.kusm.kyoto-u.ac.jp/\textasciitilde nakajima}
\thanks{Supported by the Grant-in-aid
for Scientific Research (No.13640019), JSPS}
%
\dedicatory{Dedicated to Professor Takushiro Ochiai on his sixtieth
birthday}
\begin{abstract}
We prove the Kirillov-Reshetikhin conjecture concerning certain finite
dimensional representations of a quantum affine algebra $\Ua$ when
$\widehat\g$ is an untwisted affine Lie algebra of type $ADE$. We use
$t$--analog of $q$--characters introduced by the author in an
essential way.
\end{abstract}
\maketitle
\section{Introduction}
Let $\g$ be a complex simple Lie algebra of type $ADE$, $\Lg = \g
\otimes \C[z,z^{-1}]$ be its loop algebra, and $\Ul$ be its quantum
universal enveloping algebra, or the quantum loop algebra for short.
It is a subquotient of the quantum affine algebra $\Ua$, i.e., without
central extension and degree operator.
It is customary to define $\Ul$ as an algebra over $\Q(q)$, but here
we consider $q$ as a nonzero complex number which is not a root of
unity.
Let $\alpha_i$ and $\Lambda_i$ ($i\in I$) be simple roots and
fundamental roots of $\g$, where $I$ is the index set. Let
$(a_{ij})_{i,j\in I}$ be the Cartan matrix.
By Drinfeld \cite{Dr}, Chari-Pressley \cite{CP-rep}, isomorphism
classes of irreducible finite dimension representations of $\Ul$ are
parametrized by $I$-tuples of polynomials $P = (P_i(u))_{i\in I}$ with
normalization $P_i(0) = 1$. They are called {\it Drinfeld
polynomials}. Let us denote by $L(P)$ the corresponding irreducible
representation.

For $i\in I$, $k\in\Z_{\ge 0}$, $a\in\C^*$, let
\begin{equation*}
   \PKr ik0 = \left( (\PKr ik0)_j(u) \right)_{j\in I} ;
   \qquad
   (\PKr ik0)_j(u) \defeq
   \begin{cases}
     \displaystyle\prod_{s=1}^k (1 - aq^{2s-2} u) & \text{if $j=i$},
     \\
     1 & \text{otherwise}.
   \end{cases}
\end{equation*}
Let $\Kr ik0 \defeq L(\PKr ik0)$ be the irreducible finite dimensional
representation with Drinfeld polynomial $\PKr ik0$. This particular
class of representations are called {\it Kirillov-Reshetikhin
modules}. They were introduced in \cite{KR} and have been studied
intensively by a motivation from the so-called Bethe ansatz (see
\cite[\S5.7]{KNT} and the references therein). Our first goal is to
prove the following recursion formula, called the $T$-system, which
was conjectured by Kuniba-Nakanishi-Suzuki~\cite{KNS}, and a
convergence property:
\begin{Theorem}\label{thm:main}
\textup{(1)} There exists an exact sequence
\begin{equation*}
   0 \to \bigotimes_{j:a_{ij}=-1} \Kr jk1 \to \Kr ik0 \otimes \Kr ik2
   \to \Kr i{k+1}0 \otimes \Kr i{k-1}2 \to 0
   \quad (k=1,2,\dots).
\end{equation*}
%
Moreover, the first and third terms are irreducible.
Here we have the following convention: we introduce some ordering
among $j$'s such that $a_{ij}=-1$ in order to define the tensor product
of the first term. We set $\Kr ik0$ to be the trivial representation
if $k=0$. If $\g$ is of type $A_1$, the first term of the exact sequence
is understood as the trivial representation.

\textup{(2)} The normalized $q$--character of $\Kr ik0$, considered as 
a polynomial in $A_{i,a}^{-1}$, has a limit as a formal power
series\textup:
\begin{equation*}
    \exists\lim_{k\to\infty} 
    \frac{ \qch(\Kr ik0)}{\Yq i0 \Yq i2\dots \Yq i{2k-2}}
    \in \Z[[A_{i,aq^{m}}^{-1}]]_{i\in I,m\in\Z_{>0}}.
\end{equation*}
\end{Theorem}

The definition of $\qch$ will be recalled in \secref{sec:qch}.
Taking the $\qch$ of the exact sequence, we get
\begin{equation*}
   \qch(\Kr ik0) \qch(\Kr ik2)
   = \qch(\Kr i{k+1}0) \qch(\Kr i{k-1}2)
   + \prod_{j:a_{ij}=-1} \qch(\Kr jk1)
   \quad (k=1,2,\dots).
\end{equation*}
This is the original {\it $T$-system}.
The $T$-system is a recursion formula. In particular, we can determine
all $\qch(\Kr ik0)$ inductively, once we know $\qch(\Kr i10)$ for any
$i\in I$.
These representations $\Kr i10$, i.e., Kirillov-Reshetikhin modules
with $k=1$, are called {\it l\/}--fundamental representations (or
simply `fundamental representations'). They are basic building blocks
in the representation theory of $\Ul$, since any irreducible
representations are subquotients of their tensor products. If $\g$ is
of type $A_n$, they are lifts of the representation of fundamental
representations of $\Uq$.
For general $\g$, there is a combinatorial algorithm to compute them
due to Frenkel-Mukhin~\cite{FM}. But the algorithm is too complicated
and author's computer program, written by C, did not give answers for
two fundamental representations of $E_8$ so far. For type $A_n$, $D_n$
see \cite{Na:AD}.
Even if we can compute $\qtch(\Kr i10)$ for the special node $i$ of
$E_8$, it is a polynomial consisting of 6899079264 monomials. It is
pratically impossible to compute even $\qtch(\Kr i20)$ by the
recursion.

However we can deduce a weaker information for $\Kr ik0$ as follows.
The quantum loop algebra $\Ul$ contains the finite dimensional quantum
enveloping algebra $\Uq$ as a subalgebra. Let $\operatorname{Res}$ be
the functor sending $\Ul$-modules to $\Uq$-modules by restriction. We
set $\QK ik \defeq \operatorname{Res} \Kr ik0$. It is known that it is
independent of $a$.




\thmref{thm:main} has the following specialization:
\begin{equation}\label{eq:Qsys}
    \QK ik \otimes \QK ik
   = \left(\QK i{k+1} \otimes \QK i{k-1}\right)
   \oplus \bigotimes_{j:a_{ij}=-1} \QK jk
   \quad (k=1,2,\dots)
\end{equation}
and
\begin{equation*}
    \exists\lim_{k\to\infty} e^{-k\Lambda_i}\chi(\QK ik)
    \in \Z[[e^{-\alpha_i}]]_{i\in I}.
\end{equation*}
The equation \eqref{eq:Qsys} is called {\it $Q$-system}. It also has
a recursive structure, and we do not know the initial term.
However, Hatayama-Kuniba-Okado-Takagi-Yamada \cite{HKOTY} showed that
the solution of the $Q$-system with the above convergence property,
which is a sum of characters of representations of $\g$, has the
following {\it explicit\/} expression. (See also \cite{KNT} for more
general results in this derection.)

\begin{Corollary}[Kirillov-Reshetikhin conjecture]\label{cor:KR}
Let 
\(
   \mathcal Q^{(i)}_k \defeq e^{-k\Lambda_i} \chi(\QK ik)
\)
be the normalized character of $\QK ik$.
For a sequence $\nu = (\nu_k^{(i)})_{i\in I,k\in \Z_{>0}}$ such that
all but finitely many $\nu_k^{(i)}$ are zero, we set
\(
   \mathcal Q^{\nu} \defeq \prod_{i,k}
      \left(\mathcal Q^{(i)}_k\right)^{\nu_k^{(i)}}.
\)
Then we have
\begin{gather*}
   \mathcal Q^{\nu} \prod_{\alpha\in\Delta_+} (1 - e^{-\alpha})
   = \sum_{N = (N^{(i)}_k)} \prod_{i\in I,k\in\Z_{>0}}
      \binom{P_k^{(i)}(\nu,N) + N_k^{(i)}}{N_k^{(i)}}
      e^{-kN_k^{(i)}\alpha_i},
\\
   P_k^{(i)}(\nu,N) \defeq
   \sum_{l=1}^\infty \nu_l^{(i)} \min(k,l)
     - \sum_{j\in I} a_{ij} \sum_{l\in\Z_{>0}} N_l^{(j)} \min(k,l),
\end{gather*}
where $\Delta_+$ is the set of the postive roots, and
$\tbinom{a}{b}=\Gamma(a+1)/\Gamma(a-b+1)\Gamma(b+1)$.
\end{Corollary}

\corref{cor:KR} and the equation~\eqref{eq:Qsys} were stated in
\cite{KR}, but proofs were not provided.
Then the paper became the starting point of the subject.
See \cite[5.7]{KNT} for the status of the conjecture in the last year.

For a geometric side, it gives us an explicit formula of the Euler
number (but not of Betti numbers) of arbitrary quiver varieties of
type $ADE$. The formula are slightly different from the specialization
of an explicit formula for Betti numbers, which was conjectured by
Lusztig~\cite{Lu:ferm} in connection with fermionic formula. The
difference is the definition of binomial coefficients. His ($q$-analog
of) $\tbinom{a}{b}$ is set $0$ if $b > a$. The equivalence between two
formulae are not known so far. (See \cite[Remark~1.3]{KNT} for
detail.)


For the proofs of above results, we use `$t$--analogs of
$q$--characters', introduced by the author
\cite{Na-qchar,Na-qchar-main} by using quiver varieties. They were
defined for {\it arbitrary\/} representations of $\Ul$ and there are
combinatorial algorithms to compute them for arbitrary irreducible
representations. Therefore we can, in principle, compute many things,
such as tensor product decompositions, etc. (The original
$q$--characters were introduced and studied by Knight,
Frenkel-Reshetikhin, Frenkel-Mukhin \cite{Kn,FR,FM}.) A new point in
this paper is that the algorithms are drastically simplified, when
they are applied to Kirillov-Reshetikhin modules.
In \cite[\S10]{Na-qchar-main} it was conjectured that a geometric
reason for this simplification was {\it smallness\/} of certain
natural morphisms between two graded quiver varieties. We do not prove 
this conjecture, but find that checking the smallness only for few
cases is enough to derive the $T$-system.

It is highly desirable to have a characterization of the solution of
$T$-system and its $t$-analog, similar to the result of
Kuniba-Nakanishi-Tsuboi~\cite{KNT}. Such a result should have many
applications to the representation theory of $\Ul$, according to the
importance of $t$--analogs of $q$--characters.

\subsection*{Acknowledgement}
Most of this work was done while I was visiting the Mathematical
Sciences Research Institute at Berkeley.
When I was a student of Professor Ochiai, he often talked me how
wonderful Berkeley is. Therefore it is a great pleasure for me to
dedicate this paper for his birthday to prove that~!

\section{Review on $t$--analogs of $q$--characters}\label{sec:qch}

In this section we recall the theory of $t$--analogs of
$q$--characters \cite{Na-qchar,Na-qchar-main}, which we will use
later. We give their properties, axiomized in \cite{Na-qchar-main}, so
that we do not need to go back to the original definitions.
But in order to understand a flavor of the theory, we first briefly
summarize results of \cite{FR,FM} and
\cite{Na-qaff,Na-qchar,Na-qchar-main}. This part is expository. See
\cite{Na-qchar} for more detailed survey.

Let $\bfR$ be the Grothendieck ring of the category of finite
dimensional representations of $\Ul$. The class $\{ L(P) \}_P$ of
irreducible finite dimensional representations is a basis of $\bfR$.
The $q$--character $\qch(V)$ \cite{Kn,FR,FM} is the generating
function of {\it l\/}--weight multiplicities of a representation $V$.
Here an {\it l\/}--weight space of $V$ is a simultaneous generalized
eigenspace with respect to the commutative subalgebra of $\Ul$ which
corresponds to $\mathbf U(\mathfrak h\otimes\C[z,z^{-1}])$ at $q=1$.
($\mathfrak h$ is the Cartan subalgebra of $\g$.)
The corresponding eigenvalues, which we call {\it l\/}--weights, are
parametrized by $I$-tuple of rational functions $Q/R =
(Q_i(u)/R_i(u))_{i\in I}$ with normalization $Q_i(0) = R_i(0) = 1$
\cite{FR}. We say an {\it l\/}--weight is {\it l\/}--dominant when
$Q_i(u)/R_i(u)$ are polynomials for all $i$.  These concepts are
natural analogs of the corresponding concepts without `{\it l\/}--'
for finite dimensional representations of $\g$.  That is, weight
spaces are simultaneous eigenspaces of $\mathfrak h$, weights are
their eigenvalues. Irreducible representations are parametrized by
dominant weights, and so on. There is a natural analog of the
dominance order on weights for {\it l\/}--weights. It will be denoted
by $<$.

As in the case of finite dimensional representations of $\g$, the
$q$--characters define a ring homomorphism (see \cite{FR} for the
proof)
\begin{equation*}
    \qch \colon \bfR \to \mathscr Y \defeq 
   \Z[Y_{i,a}, Y_{i,a}^{-1}]_{i\in I, a\in\C^*},
\end{equation*}
where $\mathscr Y$ is the Laurent polynomial ring of uncontably many
variables $Y_{i,a}$'s ($i\in I$, $a\in\C^*$). Moreover, this is
injective.
Then, Frenkel-Mukhin gave a characterization of the image of $\qch$,
which gives a combinatorial algorithm to compute $\qch(L(P))$ when
$L(P)$ is an {\it l\/}--fundamental representation \cite{FM}.
%

Around the same time, the author introduced another family of finite
dimensional representations $\{ M(P) \}$ of $\Ul$, also parametrized
by Drinfeld polynomial $P$ \cite{Na-qaff}. They are called {\it
standard modules}. Those $M(P)$ and irreducible representations $L(P)$
naturally define classes in the $t$--analog of the representation ring
$\bfR_t \defeq \bfR\otimes \Z[t,t^{-1}]$. We denote them by the same
symbol. There is a bar involution
$\setbox5=\hbox{A}\overline{\rule{0mm}{\ht5}\hspace*{\wd5}}$ on
$\bfR_t$. Then a main result of \cite{Na-qaff} (although stated in a
different form there) is the following characterization of the class
$L(P)$:
\begin{equation}\label{eq:can}
    \overline{L(P)} = L(P), \qquad
    L(P) \in M(P) + \sum_{Q: Q<P} t^{-1}\Z[t^{-1}] M(Q).
\end{equation}
This is very similar to definitions of Kazhdan-Lusztig bases of Hecke
algebras, canonical bases of quantum enveloping algebras for type
$ADE$, etc. As these bases, this characterization gives us an
algorithm to compute the transition matrix between two bases $\{
L(P)\}$, $\{ M(P) \}$ of $\bfR_t$, once we can compute
$\overline{M(P)}$.

In order to compute $\overline{M(P)}$, the author introduced 
the $t$--analogs of $q$--characters, which gives a $\Z[t,t^{-1}]$-linear
homomorphism
\begin{equation*}
    \qtch \colon \bfR_t \to {\mathscr Y}_t,
\end{equation*}
where
\(
   \mathscr Y_t \defeq 
   \Z[t,t^{-1},Y_{i,a}, Y_{i,a}^{-1}]_{i\in I, a\in\C^*}.
\)
This is injective again and the image has a similar characterization
as $\qch$. In particular, we have a combinatorial algorithm to compute
them for {\it l\/}--fundamental representations. Although $\qtch$ is
not a ring homomorphism, we have a simple `twist' of multiplication so
that we can compute $\qtch(M(P))$ for arbitrary $P$, which are tensor
products of {\it l\/}--fundamental representations \cite{VV1}. The bar
involution on $\qtch$ is easily defined so that $\overline{M(P)}$ can
be computed, if we know all values of $\qtch(M(P))$. Combining with
the above characterization of $L(P)$, we have a combinatorial
algorithm to compute the transition matrix, although it is practically
difficult to compute them explicitly, as we mentioned in the
introduction. This is the end of our exposition.

Now we fix notations.
We use symbols $Y_{i,a}$ following \cite{FR}. It corresponds to the
character of the fundamental weight $y_i = e^{\Lambda_i}$ in the above
analogy.
Let $\mathcal M$ be the set of monomials in
$\mathscr Y$. 
For an $I$-tuple of rational functions $Q/R = (Q_i(u)/R_i(u))_{i\in
I}$ with $Q_i(0) = R_i(0) = 1$, we set
\begin{equation*}
   e^{Q/R} \defeq
   \prod_{i\in I} \prod_{\alpha} \prod_{\beta}
     Y_{i,\alpha} Y_{i,\beta}^{-1},  
\end{equation*}
where $\alpha$ (resp.\ $\beta$) runs roots of $Q_i(1/u) = 0$
(resp.\ $R_i(1/u) = 0$), i.e.,
$Q_i(u) = \prod_\alpha ( 1 - \alpha u)$ (resp.\ $R_i(u) = \prod_\beta
(1 - \beta u)$). As a special case, an $I$-tuple of polynomials $P =
(P_i(u))_{i\in I}$ defines $e^P = e^{P/1}$.
In this way, the set $\mathcal M$ of monomials are identified with the 
set of $I$-tuple of rational functions, and the set of {\it
l\/}--dominant monomials are identified with the set of $I$-tuple of
polynomials.
A {\it monomial\/} in $\mathscr Y_t$ means a monomial only in
$Y_{i,a}^\pm$, containing no $t$'s. So they are same as monomials in
$\mathscr Y$. We consider $\Z[t,t^{-1}]$ as coefficients.

Let
\begin{equation*}
   A_{i,a} \defeq Y_{i,aq} Y_{i,aq^{-1}}
     \prod_{j:j\neq i} Y_{j,a}^{a_{ij}},
\end{equation*}
where $a_{ij}$ is the $(i,j)$-entry of the Cartan matrix. (This is
different from the original definition in \cite{FR} for general
$\g$. But it is the same when $\g$ is of type $ADE$.)

\begin{Definition}
(1) For a monomial $m\in\mathcal M$, we define $u_{i,a}(m)\in\Z$ be the
degree in $Y_{i,a}$, i.e.,
\begin{equation*}
   m = \prod_{i,a} Y_{i,a}^{u_{i,a}(m)}.
\end{equation*}

(2) A monomial $m\in\mathcal M$ is said $i$--dominant if
$u_{i,a}(m)\ge 0$ for all $a$. It is said {\it l--dominant} if
it is $i$--dominant for all $i$.

(3) Let $m, m'$ be monomials in $\mathcal M$. We say $m \le m'$ if
$m/m'$ is a monomial in $A_{i,a}^{-1}$ ($i\in I$, $a\in\C^*$).
Here a monomial in $A_{i,a}^{-1}$ means a product of nonnegative
powers of $A_{i,a}^{-1}$. It does not contain any factors
$A_{i,a}$. In such a case we define $v_{i,a}(m, m')\in\Z_{\ge 0}$ by
\begin{equation*}
   m = m' \prod_{i,a} A_{i,a}^{-v_{i,a}(m,m')}.
\end{equation*}
This is well-defined since the $q$-analog of the Cartan matrix is
invertible. We say $m < m'$ if $m\le m'$ and $m\neq m'$.

(4) For an $i$--dominant monomial $m\in\mathcal M$ we define
\begin{equation*}
   E_i(m) \defeq
    m\, \prod_a
     \sum_{r_a=0}^{u_{i,a}(m)}
     t^{r_a(u_{i,a}(m)-r_a)}
     \begin{bmatrix}
       u_{i,a}(m) \\ r_a
     \end{bmatrix}_t A_{i,aq}^{-r_a},
\end{equation*}
where
\(
\left[\begin{smallmatrix}
  n \\ r
\end{smallmatrix}\right]_t
\)
is the $t$-binomial coefficient. We call $E_i(m)$ an {\it expansion\/} 
at $m$.

(5) Let $\mathscr K_{t,i}$ be the $\Z[t,t^{-1}]$-linear
subspace of $\mathscr Y_t$ generated by $E_i(m)$ with $i$--dominant
monomials $m$. Let $\mathscr K_t$ be the intersection
$\bigcap_i \mathscr K_{t,i}$.
\end{Definition}

By the identification between monomials in $\mathcal M$ and $I$-tuples
of rational functions, we also apply the above definitions to the
latter. For example, $Q/R \le Q'/R'$ means $e^{Q/R} \le e^{Q'/R'}$.

For $m \in \mathcal M$, we define $(\tilde u_{i,a}(m))_{i\in I,
a\in\C^*}$ as the solution of
\begin{equation*}
   u_{i,a}(m) = \tilde u_{i,aq^{-1}}(m) + \tilde u_{i,aq}(m)
                  - \sum_{j: a_{ij} = -1} \tilde u_{j,a}(m).
\end{equation*}
To solve the system, we may assume that $u_{i,a}(m)=0$ unless $a$ is a
power of $q$. Then the above is a recursive system, since $q$ is not a
root of unity. So it has a unique solution such that $\tilde
u_{i,q^s}(m) = 0$ for sufficiently small $s$. Note that $\tilde
u_{i,a}(m)$ is nonzero for possibly infinitely many $a$'s, although
$u_{i,a}(m)$ is not.

Suppose that {\it l\/}--dominant monomials $m_{P^1}$, $m_{P^2}$ and
monomials $m^1\le m_{P^1}$, $m^2\le m_{P^2}$ are given. We define
an integer $d(m^1, m_{P^1}; m^2, m_{P^2})$ by
\begin{equation}\label{eq:d}
\begin{split}
   d(m^1, m_{P^1}; m^2, m_{P^2})
   & \defeq
   \sum_{i,a} \left( v_{i,aq}(m^1, m_{P^1}) u_{i,a}(m^2)
   + u_{i,aq}(m_{P^1}) v_{i,a}(m^2, m_{P^2})\right)
\\
   & = 
   \sum_{i,a} \left( u_{i,a}(m^1) v_{i,aq^{-1}}(m^2, m_{P^2})
   + v_{i,a}(m^1, m_{P^1})u_{i,aq^{-1}}(m_{P^1}) \right)
.
\end{split}
\end{equation}
We also define
\begin{equation}
\begin{split}
   &
   \tilde d(m^1, m^2) \defeq 
   - \sum_{i,a} u_{i,aq}(m^1) \tilde u_{i,a}(m^2)
,
\\
   & \ve(m^1, m^2)  \defeq \tilde d(m^1, m^2) - \tilde d(m^2, m^1)
.
\end{split}
\end{equation}
Since $u_{i,a}(m^2) = 0$ except for finitely many $a$'s, this is
well-defined. Moreover, we have
\begin{equation}\label{eq:use}
\begin{gathered}
   \ve(m^1, m^2) = 
   d(m^1, m_{P^1}; m^2, m_{P^2}) - d(m^2, m_{P^2}; m^1, m_{P^1})
   + \ve(m_{P^1}, m_{P^2}).
\end{gathered}
\end{equation}
We also write $\ve(P^1, P^2)$ instead of $\ve(m_{P^1}, m_{P^2})$.

We express the property of $\qtch$ for its slightly modified version
$\tilqtch$.

\begin{Theorem}\label{thm:ind}
\textup{(1)} 
The $\tilqtch$ of a standard module $M(P)$ has a form
\begin{equation*}
   \tilqtch(M(P)) = e^P + \sum a_m(t) m,
\end{equation*}
where the summation runs over monomials $m < e^P$.

\textup{(2)}
For each $i\in I$, $\tilqtch(M(P))$ can be expressed as a linear
combination \textup(over $\Z[t,t^{-1}]$\textup) of $E_i(m)$ with
$i$--dominant monomials $m$. Moreover, the image of $\tilqtch$ is
exactly equal to $\mathscr K_t$.

\textup{(3)}
Suppose that two $I$-tuples of polynomials $P^1 = (P^1_i)$, $P^2 =
(P^2_i)$ satisfy the following condition:
\begin{equation}
\label{eq:Z}
\begin{minipage}[m]{0.75\textwidth}
\noindent   
   $a/b\notin\{ q^n \mid n\in\Z, n \ge 2\}$ for any
   pair $a$, $b$ with $P^1_i(1/a) = 0$, $P^2_j(1/b) =
   0$ \textup($i,j\in I$\textup).
\end{minipage}
\end{equation}
Then we have
\begin{equation*}
  \tilqtch(M(P^1P^2)) =
  \sum_{m^1, m^2} t^{2d(m^1, m_{P^1}; m^2, m_{P^2})}
   a_{m^1}(t) a_{m^2}(t) m^1 m^2
,
\end{equation*}
where
\(
   \tilqtch(M(P^a)) = \sum_{m^a} a_{m^a}(t) m^a
\)
with $a=1,2$.

Moreover, properties \textup{(1),(2),(3)} uniquely determine $\tilqtch$.
\end{Theorem}

Apart from the existence problem, one can consider the above
properties (1), (2), (3) as the definition of $\tilqtch$ (an
axiomatic definition). We only use the above properties, and the
reader can safely forget the original definition.

Let us explain briefly why the properties (1), (2), (3) determine
$\tilqtch$. First consider the case $M(P)$ is an {\it
  l\/}--fundamental representation. (We have $M(P) = L(P)$ in this
case.) Then one can determine $\tilqtch(M(P))$ starting from $m_P$ and
using the property (2) inductively. The graph will never contain {\it
  l\/}--dominant monomials and will stop eventually. (The idea can be
  seen in the example below.)
For general $P$, write it as $P = P^1 P^2 P^3 \cdots$ so that each
$M(P^\alpha)$ is an {\it l\/}--fundamental representation, and
the condition~\eqref{eq:Z} is met with respect to the ordering. Then
we apply (3) successively to get
\(
   \tilqtch(M(P))
\)
from
\(
   \tilqtch(M(P^\alpha))
\)
with $\alpha=1,2,\cdots$.

We now define $\qtch(M(P))$ by
\begin{equation}\label{eq:qtch}
  \qtch(M(P)) = \sum_m t^{-d(m,m_P;m, m_P)} a_m(t) m
  \qquad \text{if\ } \tilqtch(M(P)) = \sum_m a_m(t) m.
\end{equation}

We define an involution
$\setbox5=\hbox{A}\overline{\rule{0mm}{\ht5}\hspace*{\wd5}}$ on
${\mathscr Y}_t$ by $\overline{t} = t^{-1}$,
$\overline{Y_{i,a}^\pm} = Y_{i,a}^\pm$.
One can show (see \cite[\S3]{Na-qchar-main}) that we can define an
involution
$\setbox5=\hbox{A}\overline{\rule{0mm}{\ht5}\hspace*{\wd5}}$ on
$\bfR_t$ by
\begin{equation*}
   \qtch(\overline{V}) = \overline{\qtch(V)}.
\end{equation*}
(By the second statement of \thmref{thm:ind}(2), it is enough to
check that right hand side is contained in $\mathscr K_t$, after
multiplying a power of $t$. This can be checked.)

We attach to each standard module $M(P)$ an oriented colored graph as
follows. (It is a slight modification of the graph in \cite[5.3]{FR}.)
The vertices are monomials in $\tilqtch(M_P)$. We draw a colored edge
$\xrightarrow{i,a}$ from $m_1$ to $m_2$ if $m_2 = m_1 A_{i,a}^{-1}$.
We also write the coefficients of the monomials in $\tilqtch(M(P))$.
In fact, edges are determined from monomials on vertices.
Now $E_i(m)$ in \thmref{thm:ind} gives edges with color $i$ starting
from $m$ and subsequent monomials.
Similarly we also draw a graph for an irreducible representation $L(P)$.

Let us define matrices with entries in $\Z[t,t^{-1}]$:
\begin{gather*}
  c_{PQ}(t) \defeq
  \text{the coefficient of $m_Q$ in $\qtch(M(P))$},
\\
  L_{PQ}(t) \defeq \text{the coefficient of $m_Q$ in $\qtch(L(P))$},
\\
  \text{$Z_{PQ}(t)$ by $\displaystyle M(P) = \sum_Q Z_{PQ}(t) L(Q)$}.
\end{gather*}
By geometric interpretations of these polynomials
\cite[\S8]{Na-qchar-main}, all these have positivity, i.e., they are
contained in $\Z_{\ge 0}[t,t^{-1}]$.

We define a new multiplication on $\bfR_t$ by
\begin{equation*}
    m^1 \ast m^2 \defeq
    t^{\ve(m^1, m^2)} m^1 m^2.
\end{equation*}
This is noncommutative.
One can show (see \cite[\S3]{Na-qchar-main}, where $\ast$ is denoted
by $\tilde\ast$ there) that there exists a unique multiplication
$\otimes$ on $\bfR_t$ satisfying
\begin{equation*}
    \qtch(V_1\otimes V_2)
   = \qtch(V_1)\ast \qtch(V_2).
\end{equation*}
When $t=1$, this corresponds to the tensor product on $\bfR$. So there
will be no fear of the confusion.
Varagnolo-Vasserot define the multiplation in a geometric way
\cite{VV2}. Their definition implies the positivity of structure
constants, i.e., if we define $a_{PQ}^R(t)$ by
\begin{equation*}
   L(P)\otimes L(Q) = \sum_R a_{PQ}^R(t) L(R),
\end{equation*}
it has nonnegative coefficients $a_{PQ}^R(t)\in \Z_{\ge 0}[t,t^{-1}]$.

We close this section by giving a example of computation of $\Kr 120 =
L(\PKr 120)$ for $\g=A_2$.
\begin{Example}\label{ex:graph}
Let $\g$ be of type $A_2$. We put a numbering $I = \{ 1, 2\}$.

(1) The graph of the {\it l\/}--fundametal representation $\tilqtch(\Kr
  110) = \tilqtch(M(\PKr 110))$ is the following:
\begin{equation*}
\begin{CD}
  \Yq 10 @>{1,aq}>> \Yq 12^{-1} \Yq 21
  @>{2,aq^2}>> \Yq 23^{-1}.
\end{CD}
\end{equation*}
This can be checked by \thmref{thm:ind}(1),(2).
The graph of $\tilqtch(M(\PKr 112))$ is exactly the same if we replace
$a$ by $aq^2$.

(2) By \thmref{thm:ind}(3) and \eqref{eq:qtch}, we can compute the graph of
$\qtch(M(\PKr 120))$ as
\begin{equation*}
\begin{CD}
  \Yq 10 \Yq 12 @>{1,aq}>> t^{-1}\Yq 21 @>{2,aq^2}>>
    t^{-1}\Yq 1{2} \Yq 2{3}^{-1}
\\
  @V{1,aq^3}VV @. @V{1,aq^3}VV
\\
  \Yq 10 \Yq 1{4}^{-1} \Yq 2{3} @>{1,aq}>>
  \Yq 1{2}^{-1} \Yq 1{4}^{-1} \Yq 2{1} \Yq 2{3}
  @>{2,aq^2}>> t^{-1}\Yq 1{4}^{-1}
\\
  @V{2,aq^4}VV @V{2,aq^4}VV @.
\\
  \Yq 10\Yq 2{5}^{-1} @>{1,aq}>>
  \Yq 1{2}^{-1}\Yq 2{1}\Yq 2{5}^{-1}
  @>{2,aq^2}>> \Yq 2{3}^{-1} \Yq 2{5}^{-1}
\end{CD}
\end{equation*}

(3) Finally applying the characterization \eqref{eq:can}, we find
\begin{equation*}
   t^{-\ve(\PKr 110, \PKr 112)} \Kr 110 \otimes \Kr 112 =
   M(\PKr 120) = \Kr 120 + t^{-1} L(\PKr 211)
   = \Kr 120 + t^{-1} \Kr 211.
\end{equation*}
This is a special case of the $T$-system.
The graph of $\qtch(\Kr 120)$ is obtained simply from the above graph
by erasing monomials with coefficients $t^{-1}$, i.e., $\Yq 21$, $\Yq
1{2} \Yq 2{3}^{-1}$, $\Yq 1{4}^{-1}$. The resulting graph turns out to
be given by the following rule without using \eqref{eq:can}: if one
gets an {\it l\/}--dominant monomial in the expansion at $m$, then
he/she does not add it to the graph and continue the expansion.
\end{Example}

The result of this paper is proved by checking that we have the same
structure for general $\Kr ik0$.

\section{A study of right negative monomials}

We only consider monomials $m$ whose factors are of the form
$Y_{i,aq^s}$ for a fixed $a\in\C^*$ hereafter.

\begin{Definition}
(1) For a monomial $m\neq 1$, we define
\[
   \maxexp(m) = \max \left\{ s \mid\text{there is a factor $Y_{i,aq^s}$
   appearing in $m$ for some $i$} \right\}.
\]

(2) Following \cite{FM}, we say a monomial $m$ is {\it right
negative\/} if the factors $Y_{i,aq^{\maxexp(m)}}$ have nonpositive
powers for all $i$. (It must be negative for at lease one $i$ by
definition of $\maxexp(m)$.) The product of right negative monomials
is right negative. An {\it l\/}--dominant monomial is not right
negative.

\end{Definition}

In this section, we shall prove the following.
\begin{Theorem}\label{thm:induction}
\textup{(1)} $\qtch(\Kr ik0)$ does not contain right negative
monomials other than $m_{\PKr ik0}$ corresponding to the {\it
l\/}--highest weight vector.

\textup{(2)} Let $m$ be a right negative monomial appearing 
in $\qtch(\Kr ik0)$ such that $\maxexp(m)\le 2k$. Then it is
\begin{equation*}
\begin{split}
   m & = \prod_{t=0}^{s-1} \Yq i{2t} 
    \cdot
    \prod_{t=s+1}^{k} \left( 
      \Yq i{2t}^{-1}
     \prod_{j: a_{ij}=-1} \Yq j{2t-1} \right)
\\
   & =
   \Yq i0 \cdots \Yq i{2s-2}
   \Yq i{2s+2}^{-1}\cdots \Yq i{2k}^{-1}
   \prod_{j: a_{ij}=-1} \Yq j{2s+1} \cdots \Yq j{2k-1}
\end{split}
\end{equation*}
for some $s = 0, 1, \dots, k-1$. Moreover, its coefficient is equal to 
$1$.

\textup{(3)} We have
\begin{equation*}
   t^{-\ve(\PKr ik0, \PKr i1{2k})}\Kr i{k}0 \otimes \Kr i1{2k}
   = \Kr i{k+1}0 + 
   t^{-1} L(\PKr i{k-1}0 \prod_{j:a_{ij}=-1} \PKr j1{2k-1})
   \quad (k=1,2, 3,\dots).
\end{equation*}
\end{Theorem}

\begin{Remark}
(1) We do not have
\begin{equation*}
   L(\PKr i{k-1}0\prod_{j:a_{ij}=-1} \PKr j1{2k-1})
   = \Kr i{k-1}0 \otimes\bigotimes_{j:a_{ij}=-1} \Kr j1{2k-1}
\end{equation*}
in general. A counter example can be found, e.g., $\g = A_2$,
$k=2$. Therefore \thmref{thm:induction}(3) does not give us an
inductive formula for $\qtch(\Kr ik0)$ unlike the $T$-system in
\thmref{thm:main}.

(2) It may not be easy for a reader to understand the importance of
the statements of \thmref{thm:induction}(1),(2) at first sight. But
they are pratically very useful to study tensor product decompositions
of $\Kr ik0$'s, as shown in the proof of Theorems~\ref{thm:main},
\ref{thm:induction}.
As another application of \thmref{thm:induction}(1), we can apply
Frenkel-Mukhin's algorithm \cite[\S5.5]{FM} to compute $\qch(\Kr
ik0)$. (See also \cite{Na-qchar-main} for a geometric explanation.)
This is because $\qtch(\Kr ik0)$ does not have {\it l\/}--dominant
monomial other than $m_{\PKr ik0}$. This will gives us an important
application later. See Figure~\ref{fig:graph} for a part of the graph
for $\qtch(\Kr i{k+1}0)$.
\end{Remark}

We shall prove the assertions~(1), (2) simultaneously by induction on
$k$. The statement~(3) will be proved during the induction argument,
i.e., we do not need (3) for $k$ during the proof.
The assertion (1) is true for $k=1$. (See \cite[Lemma~6.5]{FM} or
\cite[Proposition 4.13]{Na-qchar-main} for a simple geometric proof.)
The assertion (2) for $k=1$ follows from \thmref{thm:ind}(2). See also
the proof of \thmref{thm:induction}(2) below.

Let us assume the assertions~(1), (2) for $k$ and prove them for
$k+1$. We consider
\begin{equation*}
   t^{-\ve(\PKr ik0, \PKr i1{2k})}\qtch(\Kr ik0 \otimes \Kr i1{2k})
   = \Yq i0 \Yq i2 \cdots \Yq i{2k-2}\cdot \Yq i{2k} + \cdots.
\end{equation*}

\begin{Lemma}\label{lem:right}
Take monomials $m$ and $m'$ from $\qtch(\Kr ik0)$ and $\qtch(\Kr
i1{2k})$ respectively such that $m m'$ is not right negative.

\textup{(1)}
We have
\begin{equation*}
\begin{split}
   & m = \Yq i0\cdots \Yq i{2k-2}
   \quad\text{or}\quad
   \prod_{t=0}^{s-1} \Yq i{2t} 
    \cdot
    \prod_{t=s+1}^{k} \left( 
      \Yq i{2t}^{-1}
     \prod_{j: a_{ij}=-1} \Yq j{2t-1} \right),
\\
   & m' = \Yq i{2k}
\end{split}
\end{equation*}
for some $s=0,1,\dots,k-1$.

\textup{(2)} The coefficient of $mm'$ in $t^{-\ve(\PKr ik0, \PKr
i1{2k})}\qtch(\Kr ik0 \otimes \Kr i1{2k})$ is equal to $1$ if $m = \Yq
i0\cdots\linebreak[0] \Yq i{2k-2}$ and $t^{-1}$ otherwise.
\end{Lemma}

\begin{proof}
Since $m m'$ is not right negative, $m$ or $m'$ is not right
negative. Suppose that $m$ is not right negative. The induction
hypothesis~(1) for $k$ implies $m = \Yq i0 \cdots \Yq i{2k-2}$.
Moreover if $m'\neq \Yq i{2k}$, then $m'$ is right negative and
$\maxexp(m') \ge  2k+2$ by \thmref{thm:ind}(2).
Therefore $m m'$ becomes right negative, and we have a contradiction
to the assumption. Thus we have $m' = \Yq i{2k}$ in this case.

Next suppose $m'$ is not right negative. Then we have $m' = \Yq
i{2k}$. We may assume that $m$ is right negative. Since $m m'$ is not
right negative, the factor $\Yq j{\maxexp(m)}$ must be equal to
$m^{\prime-1} = \Yq i{2k}^{-1}$. By the induction hypothesis~(2) for
$k$, the monomial $m$ must be of the form of the second case of the
assertion of the lemma. Thus we have completed the proof of (1) of the
lemma.

Let us compute the coefficients of $mm'$ in 
\(
  t^{-\ve(\PKr ik0, \PKr i1{2k})} \qtch(\Kr ik0 \otimes \Kr i1{2k})
\)
when $m$, $m'$ are as in the claim. If $m m' = \Yq i0\cdots
\Yq i{2k-2}\Yq i{2k}$, then the coefficient is $1$. So we assume $m$
is the second case of the assertion of the claim. By the induction
hypothesis~(2) for $k$, the coefficient of $m$ in $\qtch(\Kr ik0)$ is
$1$. The coefficient of $m'$ in $\qtch(\Kr i1{2k})$ is also
$1$ by \thmref{thm:ind}(2). We have
\begin{equation*}
\begin{split}
   & d(m, \Yq i{2k}) = v_{i,aq^{2k+1}}(m, \Yq i0\cdots \Yq i{2k-2}) = 0,
\\
    & d(\Yq i{2k}, m) = v_{i,aq^{2k-1}}(m, \Yq i0\cdots \Yq i{2k-2}) = 1.
\end{split}
\end{equation*}
Therefore the coefficient of $m m'$ is equal to $t^{-1}$, where we
have used \eqref{eq:use}.
\end{proof}

\begin{proof}[Proof of \thmref{thm:induction}(3)]
Let us define $a_P(t)$ by
\begin{equation*}
   t^{-\ve(\PKr ik0, \PKr i1{2k})}\qtch(\Kr ik0 \otimes \Kr i1{2k})
    = \sum_P a_P(t) \qtch(L(P)),
\end{equation*}
where $P$ runs all Drinfeld polynomials. Consider the coefficient of
an {\it l\/}--dominant monomial $m_Q$ in both hand sides. Since an
{\it l\/}--dominant monomial is not right negative, \lemref{lem:right}
implies
\begin{multline*}
   m_Q = \Yq i0\cdots \Yq i{2k-2} \Yq i{2k}
       = m_{\PKr i{k+1}0}
\\
   \text{or}\quad
   m_Q = \Yq i0 \cdots \Yq i{2k-4}
   \prod_{j: a_{ij}=-1} \Yq j{2k-1}
   = m_{\PKr i{k-1}0}\prod_{j:a_{ij}=-1} m_{\PKr j1{2k-1}}.
\end{multline*}
The latter is the case $s=k-1$ of \lemref{lem:right}(1). Note that
the corresponding $m_Q$ is not {\it l\/}--dominant for $s\neq k-1$.
The coefficient of $m_Q$ is $1$ or $t^{-1}$ respectively. On the other hand,
from the right hand side, it is equal to
\begin{equation*}
   \sum_P a_P(t) L_{PQ}(t).
\end{equation*}
Remember that we have $a_P(t), L_{PQ}(t)\in \Z_{\ge 0}[t,t^{-1}]$ and
$L_{PQ}(t^{-1}) = L_{PQ}(t)$, $L_{PP}(t) = 1$. Therefore, if $Q$ is
not as above, we have $a_Q(t) = 0$. And if $Q$ is as above, we have
$a_Q(t) = 1$ or $t^{-1}$ respectively. This shows
\thmref{thm:induction}(3).
\end{proof}

\begin{proof}[Proof of \thmref{thm:induction}(1)]
Take $\qtch$ of the both hand sides of \thmref{thm:induction}(3). We have
\begin{equation*}
\begin{split}
   & \qtch(L(\PKr i{k-1}0\prod_{j:a_{ij}=-1} \PKr j1{2k-1}))
\\
  =\; & \Yq i0 \cdots \Yq i{2k-4} \prod_{j: a_{ij}=-1} \Yq j{2k-1}
\\
   & + \Yq i0 \cdots \Yq i{2k-6} \Yq i{2k-2}^{-1}
     \prod_{j: a_{ij}=-1} \Yq j{2k-3} \Yq j{2k-1}
\\
   & + \Yq i0 \cdots \Yq i{2k-8} \Yq i{2k-4}^{-1}\Yq i{2k-2}^{-1}
     \prod_{j: a_{ij}=-1} \Yq j{2k-5} \Yq j{2k-3} \Yq j{2k-1}
     + \cdots,
\end{split}
\end{equation*}
where we have used the expansion (\thmref{thm:ind}(2)). Therefore any
monomial $m$ such that
\begin{enumerate}
\item it is not right negative,
\item it appears in $\qtch$ of the left hand side of
\thmref{thm:induction}(3),
\end{enumerate}
which we just classified in \lemref{lem:right}, appears in above,
except the {\it l\/}--highest weight monomial. Thus $\qtch(\Kr
i{k+1}0)$ does not contain monomials, which are not right negative,
other than $m_{\PKr i{k+1}0}$.
\end{proof}

\begin{proof}[Proof of \thmref{thm:induction}(2)]
We may assume \thmref{thm:induction}(1) for $k+1$ now. Therefore we
can apply Frenkel-Mukhin's algorithm to compute $\qtch(\Kr i{k+1}0)$.
But we will use a self-contained argument, using only the induction
hypothesis.
We want to classify all right negative monomials $m$ with
$\maxexp(m)\le 2k+2$.

Let us consider the graph $\Gamma$ of $\qtch(\Kr i{k+1}0)$.
Since $m$ is not {\it l\/}--dominant, it must have an incoming arrow
$m'\to m$. Noticing that if $m'$ is also right negative (i.e., it is
not the {\it l\/}--highest weight monomial), then we have $\maxexp(m')
\le \maxexp(m)\le 2k+2$. We can repeat this procedure until $m^{(N)} =
(\cdots((m')')\cdots)'$ ($N$-times) becomes the {\it l\/}--highest
weight monomial $\Yq i0 \cdots \Yq i{2k}$. So we study the graph
$\Gamma$ from $\Yq i0 \cdots \Yq i{2k}$ and forget monomials $m$ with
$\maxexp(m) > 2k+2$. The first part of the graph $\Gamma$ is
\begin{equation*}
   m^{(N)} = \Yq i0 \cdots \Yq i{2k}
   \xrightarrow{i,aq^{2k+1}}
   \Yq i0 \cdots \Yq i{2k-2} \Yq i{2k+2}^{-1} 
     \prod_{j: a_{ij}=-1} \Yq j{2k+1}.
\end{equation*}
There are no other arrows from $\Yq i0 \cdots \Yq i{2k}$:
other candidates, which are of the form
\[
   \Yq i0 \cdots \Yq i{2s-2} \Yq i{2s+4} \cdots \Yq i{2k}
     \prod_{j: a_{ij}=-1} \Yq j{2s+1}
\]
are {\it l\/}--dominant, and they do not appear in $\qtch(\Kr
i{k+1}0)$. Therefore $m^{(N-1)}$ is the second monomial in the graph. 
We continue the graph. Next part is
\begin{equation*}
\begin{CD}
   {\displaystyle m^{(N-1)} = \Yq i0 \cdots \Yq i{2k-2} \Yq i{2k+2}^{-1} 
     \prod_{j: a_{ij}=-1} \Yq j{2k+1}} @>{j,aq^{2k+2}}>>
\\
   @VV{i,aq^{2k-1}}V @.
\\
   {\displaystyle \Yq i0 \cdots \Yq i{2k-4} \Yq i{2k}^{-1}\Yq i{2k+2}^{-1} 
     \prod_{j: a_{ij}=-1} \Yq j{2k-1}\Yq j{2k+1}} @..
\end{CD}
\end{equation*}
We do not consider the incoming monomial of $\xrightarrow{j,aq^{2k+2}}$
since its largest exponent $l$ is $2k+3$.
There are no other outgoing arrows from $m^{(N-1)}$,
i.e., $\xrightarrow{i,aq^{2s-1}}$ with $s=1,2,\dots,k-1$.
Such an arrow goes to a monomial, which leads to a contradiction
with the induction hypothesis~(2) for $k$, as one can easily check
by considering $\Kr ik0\otimes \Kr i1{2k}$.
We continue the graph to get all
\begin{equation*}
   \Yq i0 \cdots \Yq i{2s-2}
   \Yq i{2s+2}^{-1}\cdots \Yq i{2k+2}^{-1}
   \prod_{j: a_{ij}=-1} \Yq j{2s+1} \cdots \Yq j{2k+1}
\end{equation*}
for $s=k, k-1,\dots,0$. We cannot have an outgoing arrow
$\xrightarrow{j,aq^{2t}}$ with $t=s,\dots, k$, again by a study of
$\qtch(\Kr ik0\otimes\Kr i1{2k})$. We do not consider the incoming monomial
of $\xrightarrow{j,aq^{2k+2}}$ as above.
The last monomial is
\begin{equation*}
   \Yq i2^{-1} \cdots \Yq i{2k+2}^{-1}
   \prod_{j: a_{ij}=-1} \Yq j1 \cdots \Yq j{2k+1}.
\end{equation*}
Then an outgoing arrow goes to a monomial such that its maximal
exponent is greater than $2k+2$. Thus we obtain all the monomials
satisfying the condition, and get \thmref{thm:induction}(2) for $k+1$.
(See Figure~\ref{fig:graph} for the part of the graph for $\qtch(\Kr
ik0)$.)
\end{proof}

\begin{figure}[htbp]
\begin{center}
\leavevmode
\begin{equation*}
\begin{CD}
   \Yq i0\dots \Yq i{2k} @.
\\
   @VV{i,aq^{2k+1}}V
\\
   \Yq i0\dots \Yq i{2k-2}\Yq i{2k+2}^{-1} \prod_j \Yq j{2k+1}
   @>{j,aq^{2k+2}}>> @.
\\
   @VV{i,aq^{2k-1}}V
\\
   {\vdots}
\\
   @VV{i,aq^3}V
\\
   \Yq i0 \Yq i4^{-1} \dots \Yq i{2k+2}^{-1} \prod_j \Yq j3\dots\Yq
   j{2k+1} @>{j,aq^{2k+2}}>> @.
\\
   @VV{i,aq}V
\\
   \Yq i2^{-1} \dots \Yq i{2k+2}^{-1} \prod_j \Yq j1\dots\Yq
   j{2k+1} @>{j,aq^{2k+2}}>> @.
\end{CD}
\end{equation*}
\caption{A part of the graph for $\qtch(\Kr i{k+1}0)$}
\label{fig:graph}
\end{center}
\end{figure}

\begin{proof}[Proof of \thmref{thm:main}(2)]
We already computed a part of $\qtch(\Kr i{k+1}0)$ as
Figure~\ref{fig:graph}.
The horizontal arrows, not written here, are multiplications by
$A_{j,b}$'s for $(j,b)\neq (i,aq^{2s-1})$ ($s=1,\dots,k+1$). Compare with
the graph for $\qtch(\Kr ik2)$. If we multiply $\Yq i0$, the graph is
exactly the same as the above graph except the last line. This is
because we do not touch the factor $\Yq i0$ until the last line.
Therefore we have
\begin{equation*}
   \qtch(\Kr i{k+1}0) = \Yq i0\cdot \qtch(\Kr ik2) + 
   \Yq i0\cdots \Yq i{2k} \cdot E,
\end{equation*}
where $E$ is a polynomial in $A_{i,b}^{-1}$ of degree greater than
$k$. Here we assign $A_{i,b}^{-1}$ degree $1$ for any $i$, $b$. This
implies the convergence statement of \thmref{thm:main}(2).
\end{proof}

\section{Proof of \thmref{thm:main}(1)}

In this section, we shall complete the proof of \thmref{thm:main}. Our 
first goal is to show that $\Kr i{k+1}0 \otimes \Kr i{k-1}2$ is
irreducible. We need the following.
\begin{Lemma}\label{lem:semismall}
Let $P = \PKr ik0\PKr ik2 = \PKr i{k+1}0 \PKr i{k-1}2$.
Then $\qtch(L(P))$ contains all
monomials
\begin{equation*}
\begin{aligned}[c]
   & \prod_{t=0}^{s-1} \Yq i{2t} \Yq i{2t+2} \cdot
     \prod_{t=s}^{k-1} \Yq j{2t+1}
\\   
   =\; & \Yq i0 \Yq i2^2 \cdots \Yq i{2s-2}^2 \Yq i{2s}
   \prod_{j:a_{ij}=-1} \Yq j{2s+1}\cdots \Yq j{2k-1}
\end{aligned}
\quad (s=1,2,\dots,k-1)
\end{equation*}
with coefficients $1$.
\end{Lemma}

\begin{proof}
Let us consider the standard module $M(P)$ with Drinfeld polynomial
$P$. We consider first several terms of $\qtch(M(P))$:
\begin{multline*}
   \qtch(M(P))
   = \Yq i0\Yq i2^2 \cdots \Yq i{2k-2}^2 \Yq i{2k}
\\
   + \sum_{s=1}^{k-1} (1 + t^2)^{k-s} t^{2(s-k)}
     \Yq i0\Yq i2^2 \cdots \Yq i{2s-2}^2 \Yq i{2s}
     \prod_{j: a_{ij} = -1} \Yq j{2s+1} \cdots \Yq j{2k-1}
\\
   + \cdots.
\end{multline*}
Note that 
\(
   (1 + t^2)^{k-s} t^{2(s-k)} \in \Z[t^{-1}]
\)
and the constant term is equal to $1$. Then the characterization
\eqref{eq:can} gives us
\begin{multline*}
   \qtch(L(P))
   = \Yq i0\Yq i2^2 \cdots \Yq i{2k-2}^2 \Yq i{2k}
\\
   + \sum_{s=1}^{k-1}
     \Yq i0\Yq i2^2 \cdots \Yq i{2s-2}^2 \Yq i{2s}
     \prod_{j: a_{ij} = -1} \Yq j{2s+1} \cdots \Yq j{2k-1}
\\
   + \cdots.
\end{multline*}
The procedure to define $L(P)$ from $M(P)$ is recursive. So we do not
need to cosider the further part.
\end{proof}

\begin{Remark}
During the proof, we encountered a simplification of the
characterization \eqref{eq:can}, i.e., simply taking constant terms.
This is exactly the phenomena appearing when we study the
decomposition theorem for a {\it semi-small\/} morphism $\pi\colon
X\to Y$. (See \cite[\S7]{Na-qchar}, \cite[9.3]{Na-qchar-main} for
similar arguments.)
  
It seems natural to expect that $M(P)$ is {\it semi-small\/} in the
sense of \cite[10.1]{Na-qchar-main}, i.e., $c_{QR}(t)\in \Z[t^{-1}]$
for all $Q,R\le P$, although we have checked the assertion for only
first several terms.
\end{Remark}

\begin{proof}[Proof of the irreducibility of $\Kr i{k+1}0 \otimes \Kr
i{k-1}2$]
We consider 
\[
   t^{-\ve(\PKr i{k+1}0, \PKr i{k-1}2)}
     \qtch(\Kr i{k+1}0 \otimes \Kr i{k-1}2).
\]
Our first task is to classify all {\it l\/}--dominant monomials
appearing in it.
Take monomials $m$ and $m'$ from $\qtch(\Kr i{k+1}0)$ and $\qtch(\Kr
i{k-1}2)$ respectively such that $m m'$ is not right negative. If $m$
is right negative, then $m'$ is not right negative, hence
we have $m' = \Yq i2\cdots \Yq i{2k-2}$ by
\thmref{thm:induction}(1). But then $m m'$ must be right negative as
we can see from \thmref{thm:induction}(2), applied to $m$. Therefore
$m$ is not right negative and we have $m = \Yq i0\cdots \Yq i{2k}$ by
\thmref{thm:induction}(1).
Since $mm'$ is not right negative, the factor $\Yq j{\maxexp(m')}$
must be one of $\Yq i0^{-1}$, \dots, $\Yq i{2k}^{-1}$. (In fact, $\Yq
i0^{-1}$ is impossible.) By \thmref{thm:induction}(2) we have
\begin{equation*}
\begin{split}
    & m m'
\\
   =\;& \Yq i0\cdots \Yq i{2k} \cdot
     \Yq i2 \cdots \Yq i{2s-2} \Yq i{2s+2}^{-1}\cdots \Yq i{2k}^{-1}
     \prod_{j: a_{ij} = -1} \Yq j{2s+1} \cdots \Yq j{2k-1}
\\
   =\; &  \Yq i0\Yq i2^2 \cdots \Yq i{2s-2}^2 \Yq i{2s}
     \prod_{j: a_{ij} = -1} \Yq j{2s+1} \cdots \Yq j{2k-1}
\end{split}
\end{equation*}
for some $s=1,2,\dots,k-1$. This is the classification. Again by
\thmref{thm:induction}(2), the coefficients of above $m'$ in
$\qtch(\Kr i{k-1}2)$ is $1$. A direct computation shows
\begin{equation*}
   d(m, m_{\PKr i{k+1}0}; m', m_{\PKr i{k-1}2})
   = d(m', m_{\PKr i{k-1}2}; m, m_{\PKr i{k+1}0}) = 0.
\end{equation*}
Therefore the coefficient of $mm'$ in
$t^{-\ve(\PKr i{k+1}0, \PKr i{k-1}2)}\qtch(\Kr i{k+1}0\otimes\Kr i{k-1}2)$ is also $1$.

Let us define $a_Q(t)$ by
\begin{equation}\label{eq:dec}
   \qtch(\Kr i{k+1}0 \otimes \Kr i{k-1}2) = \sum_Q a_Q(t) \qtch(L(Q)),
\end{equation}
where $Q$ runs all Drinfeld polynomials. We have $a_P(t) = 1$. Let us
consider $Q\neq P$. If $m_Q$ is not one of the above classified
monomials, then the coefficient in the left hand side is $0$. Since
$L_{QQ}(t) = 1$, the positivity of coefficients implies $a_Q(t) = 0$
for such $Q$. Next we consider $m_Q$ as above. We have
\begin{equation*}
   1 = \sum_{Q'} a_{Q'}(t) L_{Q'Q}(t).
\end{equation*}
But \lemref{lem:semismall} means $L_{PQ}(t) = 1$. Therefore the
positivity of $L_{Q'Q}(t)$ implies that $a_{Q'} L_{Q'Q}(t) = 0$ for
$Q'\neq P$. Setting $Q' = Q$, we get $a_Q(t) = 0$ since $L_{QQ}(t) =
1$. This means that the right hand side of \eqref{eq:dec} is a single
term $L(P)$, that is $\Kr i{k+1}0 \otimes \Kr i{k-1}2$ is irreducible.
\end{proof}

\begin{Remark}
  Chari gave a sufficient condition for the irreducibility of the tensor
  product of Kirillov-Reshetikhin modules \cite{Ch}. But the above
  case is not covered by her result.
\end{Remark}



\begin{proof}[Completion of the proof of \thmref{thm:main}]
First of all, it is easy to show the irreducibility of the tensor
product $\bigotimes_{j:a_{ij}=-1} \Kr jk1$: Using
\thmref{thm:induction}(1),(2), one can prove that its $\qtch$ has no
{\it l\/}--dominant monomials other than the {\it l\/}--highest weight
monomial. This case can be also deduced from Chari's sufficient condition
\cite{Ch}.

The rest of the proof goes on a similar line as the above proof of the
irreducibility of $\Kr i{k+1}0 \otimes \Kr i{k-1}2$. Let us consider
\(
   t^{-\ve(\PKr ik0, \PKr ik2)} \qtch(\Kr ik0 \otimes \Kr ik2).
\)
Take monomials $m$ and $m'$ from $\qtch(\Kr ik0)$ and $\qtch(\Kr ik2)$
respectively such that $m m'$ is not right negative. By an argument as
above, we have
\begin{equation*}
    m m'
   = \Yq i0\Yq i2^2 \cdots \Yq i{2s-2}^2 \Yq i{2s}
     \prod_{j: a_{ij} = -1} \Yq j{2s+1} \cdots \Yq j{2k-1}
\end{equation*}
for some $s=0,\dots, k-1$. 
(We have $m' = \Yq i2\cdots \Yq i{2k}$.)
Note that we have one extra monomial when $s=0$, i.e.,
\begin{equation*}
    \prod_{j: a_{ij} = -1} \Yq j1 \cdots \Yq j{2k-1}.
\end{equation*}
We have
\begin{equation*}
\begin{split}
 & d(m, m_{\PKr i{k}0}; m', m_{\PKr i{k}2})
   - d(m', m_{\PKr i{k}2}; m, m_{\PKr i{k}0})
\\
  =\; & 
   v_{i,aq^{2k+1}}(m, m_{\PKr i{k}0}) - v_{i,aq}(m, m_{\PKr i{k}0})
   = 
   \begin{cases}
     -1 & \text{if $s=0$},
     \\
     0 & \text{otherwise}.
   \end{cases}
\end{split}
\end{equation*}
Therefore the coefficient of $mm'$ in
\( 
   t^{-\ve(\PKr ik0, \PKr ik2)}\qtch(\Kr ik0 \otimes\Kr ik2)
\)
is $t^{-1}$ ($s=0$) or $1$ ($s\neq 0$). Considering the decomposition of
$\qtch(\Kr ik0 \otimes\Kr ik2)$ as above, we get
\begin{multline*}
   t^{-\ve(\PKr ik0, \PKr ik2)}\qtch(\Kr ik0 \otimes \Kr ik2) 
\\
   = t^{-\ve(\PKr i{k+1}0, \PKr i{k-1}2)}
           \qtch(\Kr i{k+1}0 \otimes \Kr i{k-1}2)
   + t^{-1-N} 
   \qtch(\bigotimes_{j:a_{ij}=-1}
   \Kr jk1),
\end{multline*}
where $N=\sum_{a < b} \ve(\PKr {j_a}k1, \PKr {j_b}k1)$. Here
$j_1, j_2, \dots$ is the ordering used for the definition of
the tensor product. If we set $t=1$, this is nothing but the $T$-system.

By \cite{Ch}, $\Kr ik0 \otimes \Kr ik2$ is cyclic: it is generated by
the tensor product of {\it l\/}--highest weight vectors. Its simple
quotient must be isomorphic to $\Kr i{k+1}0 \otimes \Kr i{k-1}2$ since 
they have the same Drinfeld polynomials. This shows the existence of
the exact sequence.
\end{proof}


\begin{thebibliography}{99}
\bibitem{Ch} V.~Chari,
{\it Braid group actions and tensor products},
preprint, arXiv:math.QA/0106241.

\bibitem{CP-rep} V.~Chari and A.~Pressley, 
\emph{Quantum affine algebras and their representations},
in {\it Representations of groups (Banff, AB, 1994)},
Amer. Math. Soc., Providence, RI, 1995, pp.~59--78.

\bibitem{Dr} V.G.~Drinfel'd,
{\it A new realization of Yangians and quantized affine algebras},
Soviet math. Dokl. {\bf 32} (1988), 212--216.

\bibitem{FR1} E.~Frenkel and N.~Reshetikhin,
\emph{Quantum affine algebras and deformations of the Virasoro and $\mathscr {W}$-algebras}, Comm. Math. Phys. {\bf 178} (1996), 237--264.

\bibitem{FR} \bysame,
\emph{The $q$--characters of representations of quantum affine algebras
  and deformations of $\mathscr {W}$--algebras},
in {\it Recent developments in quantum affine algebras and related topics
(Raleigh, NC, 1998)},
Contemp. Math., {\bf 248}, Amer. Math. Soc., Providence, RI, 1999, 163--205.

\bibitem{FM} E.~Frenkel and E.~Mukhin,
\emph{{Combinatorics of $q$--characters of
  finite-dimensional representations of quantum affine algebras}},
Comm. Math. Phys. {\bf 216} (2001), 23--57.

\bibitem{HKOTY} G.~Hatayama, A.~Kuniba, M.~Okado, T.~Takagi and
  Y.~Yamada,
{\it Remarks on fermionic formula},
in {\it Recent developments in quantum affine algebras and related topics
(Raleigh, NC, 1998)}, Contemp. Math. {\bf 248} (1999), 243--291.

\bibitem{KR} A.N.~Kirillov and N.~Reshetikhin,
{\it Representation of Yangians and multiplicity of occurrence of the
irreducible components of the tensor product of representations of
simple Lie algebras},
J. Sov. Math. {\bf 52} (1990), 3156--3164.

\bibitem{Kn} H.~Knight, \emph{Spectra of tensor products of
    finite-dimensional representations of Yangians},
J. Algebra {\bf 174} (1995), 187--196.

\bibitem{KNS} A.~Kuniba, T.~Nakanishi and J.~Suzuki,
{\it Functional relations in solvable lattice models: I. Functional
relations and representation theory},
Int. J. Mod. Phys. A {\bf 9} (1994), 5215--5266.

\bibitem{KNT} A.~Kuniba, T.~Nakanishi and Z.~Tsuboi,
{\it The canonical solutions of the $Q$-systems and the
Kirillov-Reshetikhin conjecture},
preprint, arXiv:math.QA/0105145.


\bibitem{Lu:ferm} G.~Lusztig,
\emph{{{F}erminonic form and {B}etti numbers}}, preprint,
arXiv:math.QA/0005010.




\bibitem{Na-qaff} H.~Nakajima,
\emph{{Quiver varieties and finite dimensional representations
  of quantum affine algebras}}, J. Amer. Math. Soc.,
\textbf{14} (2001), 145--238.

\bibitem{Na-qchar} \bysame,
\emph{{$t$--analogue of the $q$--characters of finite
dimensional representations of quantum affine algebras}},
in ``Physics and Combinatorics'', Proceedings of the Nagoya 2000
International Workshop, World Scientific, 2001, 195--218.


\bibitem{Na-qchar-main}
\bysame,
{\it Quiver varieties and $t$--analogs of $q$--characters of quantum
affine algebras},
preprint, arXiv:math.QA/0105173.

\bibitem{Na:AD} \bysame,
{\it $t$--analogs of $q$--characters of quantum affine algebras of type
$A_n$, $D_n$},
to appear.

\bibitem{VV1} M.~Varagnolo and E.~Vasserot,
{\it Standard modules of quantum affine algebras},
preprint, arXiv:math.QA/0006084, to appear in Duke Math. J.

\bibitem{VV2} \bysame,
{\it Perverse sheaves and quantum Grothendieck rings},
preprint, arXiv:math.QA/0103182.

\end{thebibliography}
\end{document}